\documentclass{article}

\usepackage{graphicx}
\usepackage{amsthm,xcolor}
\usepackage{amsfonts}
\usepackage{amsmath}
\usepackage{amscd}
\usepackage{amssymb}
\usepackage{alltt}
\usepackage{url}
\usepackage{ellipsis}
\usepackage{lineno}
\usepackage{setspace}
\usepackage[numbers]{natbib}

\theoremstyle{thmstyleone}%
\theoremstyle{thmstyletwo}%

\theoremstyle{thmstylethree}%

\begin{document}

\title{A review of Alfred North Whitehead's ``Introduction to Mathematics''}

\author{Thomas Hales\\University of Pittsburgh}   
\date{}
\maketitle


\maketitle


\section{An Introduction to Mathematics}

In 1911, Alfred North Whitehead published a short book
\emph{Introduction to Mathematics} (IM) intended for students wanting
an explanation of the fundamental ideas of mathematics~\cite{IM}.

Whitehead's IM has enduring value because it was written not long
after he and Bertrand Russell published their monumental three-volume
work \emph{Principia Mathematica} (PM) -- a publication of immense
historical significance for mathematics~\cite{PM}.  IM sheds light on
Whitehead's view of mathematics at that time.  For Russell, the PM was
the fulfillment of the logicist thesis: the claim that all of
mathematics is logic and nothing but logic.  However, Russell might
have claimed more than PM could deliver, because the PM included the
axiom of choice and the axiom of infinity, which are not part of
logic.  The primary components of their system were predicate logic,
relations, and types.

Russell and Whitehead both turned away from mathematics after becoming
exhausted by the ten years they spent writing the PM.  Russell wrote
in his autobiography,
``My intellect never quite recovered from the strain. I have been ever
since definitely less capable of dealing with difficult abstractions
than I was before. This is part, though by no means the whole, of the
reason for the change in the nature of my work.''  In 1910, Whitehead
resigned his Senior Lectureship at Trinity College, Cambridge, and
moved to London, where he was unemployed for a year~\cite{Lowe}. His interests
shifted from mathematics to education, and IM was his first book after
PM.  IM is a thus a retrospective and transitional work, written for
the masses by a man recovering from a long intellectual ordeal.


The book IM was published in the series ``The Home University Library
of Modern Knowledge,'' which published concise introductions on a
broad range of topics for a general audience.  We might compare this
series with similar book series published today, such as the ``Very
Short Introduction,'' ``For Dummies,'' and the French ``Que sais-je?''
Bertrand Russell also published two books in the Home University
series, including his \emph{The Problems of Philosophy}~\cite{Russell}.  (As an
aside, Russell's \emph{Problems} was the topic of one of my college entrance
essays when I was a young man.)

Although a polymath, Whitehead was a mathematical physicist above all,
and his book IM contains a mixture of math and physics.  The closest
parallel to Whitehead today might be Roger Penrose, whose books cover
a broad range of math, physics, and philosophy~\cite{Penrose}.

IM tells the famous story of Archimedes'
eureka moment of devising a test of the purity of the gold in the
king's crown by submerging it in water to measure its specific
gravity.  According to Whitehead, ``This day, if we knew which it was,
ought to be celebrated as the birthday of mathematical physics.''  IM
was written in the early days of quantum physics and relativity, but
Whitehead does not mention the new physics. Instead, IM recounts the
development of classical physics, including the fundamental forces of
gravity, electricity, and magnetism, and Newton's laws of motion.

Vectors and linear algebra were not nearly so well known by math and
engineering students in 1910 as they are today.  Whitehead's IM helped
to popularize vectors as a tool of applied mathematics, especially in
his presentation of topics such as the parallelogram rule for the
addition of forces.  Whitehead's earlier treatise on universal algebra
contains several chapters on vectors, matrices, and forces.  The
treatise developed Cayley's theory of matrices and the properties of
vectors as conceived by Grassmann and Gibbs.  According to Whitehead
in IM, ``The idea of the `vector', that is, of a directed magnitude,
is the root-idea of physical science.\,\ldots\ Thus, when in
analytical geometry the ideas of the `origin', of `co-ordinates', and
of `vectors' are introduced, we are studying the abstracts conceptions
which correspond to the fundamental facts of the physical world.''
Vectors are used to illustrate the main object of the book IM, which
is to show students why mathematics ``is necessarily the foundation of
exact thought as applied to natural phenomena.''



In Whitehead's technical writings, his meaning is sometimes obscure,
because of his insistence on generality.  However, in his books
written for a large audience such as IM, there are many memorable
sentences.  One of the most widely quoted sentences from IM foretold
the rise of automation: ``Civilization advances by extending the
number of important operations which we can perform without thinking
about them.'' Another passage foreshadows the concepts of entropy and
compression in information theory: symbolism ``should be
concise,\ldots\ Now we cannot place symbols more concisely than by
placing them in immediate juxtaposition.  In a good symbolism
therefore, the juxtaposition of important symbols should have an
important meaning.''  Whitehead saw the quest for the permanent and
the general as fundamental, ``To see what is general in what is
particular and what is permanent in what is transitory is the aim of
scientific thought.\ldots\ This possibility of disentangling the most
complex evanescent circumstances into various examples of permanent
laws is the controlling idea of modern thought.''


In IM, Whitehead tells the well-known story that Alexander the Great
impatiently asked for his geometry tutor to make the proofs shorter.
His tutor Menaechmus made the famous reply, ``In the country there are
private and even royal roads, but in geometry there is only one road
for all.''  Whitehead disagreed, ``But if Menaechmus thought that his
proofs could not be shorteneed, he was grievously
mistaken,\ldots\ Nothing illustrates better the gain in power which is
obtained by the introduction of relevant ideas into a science than to
observe the progressive shortening of proofs which accompanies the
growth of richness in idea.\ldots\ There are royal roads in science;
but those who first tread them are men of genius and not kings.''

\section{Some and Any}

Whitehead wrote in IM, ``Mathematics as a science commenced when first
someone, probably a Greek, proved propositions about \emph{any} things
or about \emph{some} things, without specification of definite
particular things.  When Whitehead mentions the words \emph{any} and
\emph{some}, he is referring to what is now called predicate logic: 
the logic of the quantifiers \emph{for all} $(\forall)$ and
\emph{there exists} $(\exists)$.  IM thus places proofs in predicate
logic as the mythical starting point of mathematics.

This statement in IM is particularly remarkable, because Whitehead
himself was slow to understand the significance of predicate logic.
Russell learned of symbolic notation for quantifiers from Peano and
Frege, who had a major influence on Russell starting in 1900, when
Russell met Peano in Paris at the International Congress of
Mathematicians.  It was Peano who introduced the symbol $\exists$ for
the existential quantifier, which was adopted in PM.  Russell and
Whitehead used $(x)$ to mean ``for all $x$,'' and it was only later
that Hilbert introduced the symbol $(\forall)$.


Surprisingly, it appears that Whitehead failed to grasp the
significance of symbolic notation for quantifiers when he wrote
\emph{Treatise on Universal Algebra} (UA), which was published in
1898~\cite{UA}.  UA developed algebra and geometry from axioms with
precise definitions, very much in the spirit of Bourbaki a few decades
later.  The treatise is divided into books, and Book II is \emph{The
  Algebra of Symbolic Logic}.  Within that book, Whitehead has an
entire chapter on \emph{existential expressions}.  From his citations
in Book II, it is known that Whitehead was well aware of the research on
logic by Charles Saunders Pierce and his students.  Pierce had
developed notation for quantifiers by
1883~\cite{sep-logic-firstorder-emergence}.  Whitehead also cited
Schr\"oder's work, which introduced Pierce's ideas to Europe.
However, it seems that Whitehead truly failed to grasp the
significance of authors he cited in UA until after Russell's encounter
with Peano in 1900.

Whitehead endorsed the logicist thesis by writing, ``A notable
discovery has been made by the joint and partially independent work of
three men: Frege, Peano and Bertrand Russell --a German, an Italian
and an Englishman--\ldots\ of the generalized conception of the
variable and of its essential presence in all mathematical reasoning.
This discovery empties mathematics of everything but its logic''~\cite{PhiloMath}.
In this context, when Whitehead speaks of the variable, he means the
bound variable in logic.  He explained explicitly in IM, ``The
idea of the undetermined `variable' as occurring in the use of `some'
or `any' is the really important one in mathematics.''  Whitehead
admitted in IM that a full understanding of quantifiers had come only
recently, ``It was not till within the last few years that it has been
realized how fundamental \emph{any} and \emph{some} are to the very
nature of mathematics.''



\section{Relations and Axioms}

As mentioned above, one of the main components of the logic of PM was
a theory of relations or predicates.  In logic, a predicate is
interpreted as a relation of zero, one, or more variables.  A function
is a special kind of relation.  Over a period of many years, Whitehead
developed various axiomatic systems of geometry (and by extension, the
natural world). Each axiomatic system can be described as a collection
of entities that are joined by relations that are governed by axioms.  If we
try to summarize Whitehead's life endeavor in a single sentence, we
might turn to a sentence in IM: ``The events of this ever-shifting
world are but examples of a few general connexions or relations called
laws.''  

Whitehead was part of the axiomatization movement in mathematics that
became widespread through Hilbert's influence.  In the ancient and
medieval world, the axiomatic reasoning of Euclid's \emph{Elements}
served as the paradigm.  By the nineteenth century, logic and
mathematical rigor had advanced far beyond Euclid's ancient
standards. Mathematicians such as Pasch and Peano started to redo
Euclid according to more rigorous standards.  In 1899, Hilbert
published \emph{Grundlagen der Geometrie (Foundations of Geometry)},
which brought Euclid into the twentieth century with a fresh set of
definitions and axioms~\cite{Hilbert}.

According to Eves, ``By developing a postulate set for Euclidean
geometry that does not depart too greatly in spirit from Euclid's own,
and by employing a minimum of symbolism, Hilbert succeeded in
convincing mathematicians to a far greater extent than had Pasch and
Peano, of the purely hypothetico-deductive nature of geometry. But the
influence of Hilbert's work went far beyond this, for, backed by the
author's great mathematical authority, it firmly implanted the
postulational method, not only in the field of geometry, but also in
essentially every other branch of mathematics. The stimulus to the
development of the foundations of mathematics provided by Hilbert's
little book is difficult to overestimate''~\cite{Eves}.

Eventually, general foundations for mathematics were provided by the
type theory of PM (and also, independently by the Zermelo axioms of
set theory).  These foundational systems have had an enormous
influence on the growth of mathematics.  Before then, each branch of
mathematics was separately axiomatized.  To axiomatize a domain of
mathematics means to systematize it by deciding which relations are to
be considered as fundamental and by enunciating a set of axioms from
which the entire theory can be deduced.

During the decade that he and Russell devoted to the PM, as a
participant in the axiomatization movement, Whitehead wrote books on
the axioms of geometry and physics~\cite{ADG}\cite{APG}.
Whitehead's aim in one publication during that period was to
generalize the axioms of static geometry to give a system of axioms
that describe the dynamical laws of physics~\cite{MC}.  Each system in
that book gives relations on a collection of undefined entities, subject to a
list of axioms.

Even his earlier treatise on universal algebra was a project in
axiomatization of abstract algebra.  For example, the treatise
contains the axioms for a general nonassociative semiring.  Theorems
were derived from the definitions and axioms in a style of abstract
algebra that was brought into the mainstream of mathematics decades
later by B. L. van der Waerden in his book \emph{Moderne Algebra}~\cite{vdW}.






\section{Whitehead's Three Tests of Accuracy}



How is an axiomatic mathematical theory to be judged?  Whitehead suggested three
tests of accuracy:
\emph{``The tests of accuracy are logical coherence, adequacy, and
  exemplification}''\cite{RM}.%
\footnote{Whitehead proposed these tests to 
evaluate metaphysical theories, but here the tests are applied to
mathematical theories.}  
These words were written in 1926, several years before G\"odel
published his incompletness theorems, but Whitehead's three tests of
accuracy might be understood more fully in a post-G\"odelian context.
(Recall that G\"odel's first version of the incompleteness theorem was
written about the mathematical system developed by Russell and
Whitehead.) The first test of accuracy is logical coherence, which 
might be taken to mean syntactic consistency.  According to Whitehead,
only this first test is required to establish validity of
reasoning~\cite[p.~vi]{UA}.  Adequacy is a pragmatic test.  In a
mathematical system, adequacy might mean that it contains sufficient
arithmetic to carry out the constructions of G\"odel's incompleteness
theorem.  In a general metaphysical system, adequacy might mean
adequacy for the purpose of sustaining an advanced scientific
civilization.  Exemplification means that the theory admits a model
(in the sense of model theory).  The first and third tests, logical
coherence and exemplification, are related by Henkin's theorem (1949):
if a theory is syntactically consistent, then it has a model.


As suggested by this short review, Whitehead developed themes in IM
that have enduring value.


\bibliographystyle{plainnat}
\bibliography{all}

\begin{thebibliography}{15}
\providecommand{\natexlab}[1]{#1}
\providecommand{\url}[1]{\texttt{#1}}
\expandafter\ifx\csname urlstyle\endcsname\relax
  \providecommand{\doi}[1]{doi: #1}\else
  \providecommand{\doi}{doi: \begingroup \urlstyle{rm}\Url}\fi

\bibitem[Eves(1963)]{Eves}
Howard Eves.
\newblock \emph{A Survey of Geometry}, volume~I.
\newblock Allyn and Bacon, 1963.
\newblock quoted in
  \url{https://en.wikipedia.org/wiki/Foundations_of_geometry}.

\bibitem[Ewald(2019)]{sep-logic-firstorder-emergence}
William Ewald.
\newblock {The Emergence of First-Order Logic}.
\newblock In Edward~N. Zalta, editor, \emph{The {Stanford} Encyclopedia of
  Philosophy}. Metaphysics Research Lab, Stanford University, {S}pring 2019
  edition, 2019.

\bibitem[Hilbert(1899)]{Hilbert}
David Hilbert.
\newblock \emph{Grundlagen der Geometrie}.
\newblock Springer-Verlag, 1899.

\bibitem[Lowe(1985, 1990)]{Lowe}
Victor Lowe.
\newblock \emph{Alfred North Whitehead, The Man and His Work}, volume I,II.
\newblock Johns Hopkins, 1985, 1990.

\bibitem[Penrose(2004)]{Penrose}
Roger Penrose.
\newblock \emph{The Road to Reality}.
\newblock Jonathan Cape, 2004.

\bibitem[Russell(1912)]{Russell}
Bertrand Russell.
\newblock \emph{The Problems of Philosophy}.
\newblock Williams and Norgate, 1912.

\bibitem[Russell and Whitehead(1910, 1912, 1913)]{PM}
Bertrand Russell and Alfred~N Whitehead.
\newblock \emph{Principia Mathematica}.
\newblock Cambridge University Press, 1910, 1912, 1913.

\bibitem[van~der Waerden(1930)]{vdW}
B.~L. van~der Waerden.
\newblock \emph{Moderne Algebra}.
\newblock Springer-Verlag, 1930.

\bibitem[Whitehead(1898)]{UA}
Alfred~N Whitehead.
\newblock \emph{A Treatise on Universal Algebra}.
\newblock Cambridge University Press, 1898.

\bibitem[Whitehead(1906{\natexlab{a}})]{APG}
Alfred~N Whitehead.
\newblock \emph{The Axioms of Projective Geometry}.
\newblock Cambridge Tracts in Mathematics and Mathematical Physics,
  1906{\natexlab{a}}.

\bibitem[Whitehead(1906{\natexlab{b}})]{MC}
Alfred~N Whitehead.
\newblock On mathematical concepts of the material world.
\newblock \emph{Philosophical transactions of the Royal Society of London, A},
  205:\penalty0 465--525, 1906{\natexlab{b}}.

\bibitem[Whitehead(1907)]{ADG}
Alfred~N Whitehead.
\newblock \emph{The Axioms of Descriptive Geometry}.
\newblock Cambridge Tracts in Mathematics and Mathematical Physics, 1907.

\bibitem[Whitehead(1910)]{PhiloMath}
Alfred~N Whitehead.
\newblock The philosophy of mathematics.
\newblock \emph{Science progress in the twentieth century}, 5:\penalty0
  234--239, 1910.

\bibitem[Whitehead(1911)]{IM}
Alfred~N Whitehead.
\newblock \emph{An Introduction to Mathematics}.
\newblock Williams and Norgate, 1911.

\bibitem[Whitehead(1926)]{RM}
Alfred~N Whitehead.
\newblock \emph{Religion in the Making}.
\newblock Macmillan, 1926.

\end{thebibliography}

\end{document}